\documentclass[12pt,oneside,english]{amsart}
\usepackage[T1]{fontenc}
\usepackage[latin9]{inputenc}
\usepackage{amstext}
\usepackage{amsthm}
\usepackage{amssymb}
\usepackage{graphicx}

\makeatletter
\numberwithin{equation}{section}
\numberwithin{figure}{section}
\theoremstyle{plain}
\newtheorem{thm}{\protect\theoremname}
\theoremstyle{remark}
\newtheorem{rem}[thm]{\protect\remarkname}
\theoremstyle{definition}
\newtheorem{defn}[thm]{\protect\definitionname}
\theoremstyle{plain}
\newtheorem{prop}[thm]{\protect\propositionname}

\usepackage{pdflscape}
\usepackage{babel}

\providecommand{\definitionname}{Definition}
\providecommand{\remarkname}{Remark}
\providecommand{\theoremname}{Theorem}

\providecommand{\propositionname}{Proposition}

\makeatother

\usepackage{babel}
\providecommand{\definitionname}{Definition}
\providecommand{\propositionname}{Proposition}
\providecommand{\remarkname}{Remark}
\providecommand{\theoremname}{Theorem}

\begin{document}
\title{Rao distances and Conformal Mapping}

\maketitle
\begin{center}
\author{Arni S.R. Srinivasa Rao{*},}
\begin{center}
Laboratory for Theory and Mathematical Modeling, 
\par\end{center}

\begin{center}
Medical College of Georgia, 
\par\end{center}

\begin{center}
Department of Mathematics, 
\par\end{center}

\begin{center}
Augusta University, Georgia, USA 
\par\end{center}

\begin{center}
email: arni.rao2020@gmail.com 
\par\end{center}

{*}Corresponding author 
\begin{center}
\vspace{0.3cm}
 
\par\end{center}
\author{Steven G. Krantz,}
\begin{center}
Department of Mathematics, 
\par\end{center}

\begin{center}
Washington University in St. Louis, Missouri, USA 
\par\end{center}

\begin{center}
email: sgkrantz@gmail.com 
\par\end{center}

\begin{center}
\vspace{0.3cm}
\par\end{center}

\end{center}
\begin{abstract}
In this article, we have described the Rao distance (due to C.R. Rao)
and ideas of conformal mappings on 3D objects with angle preservations.
Three propositions help us to construct distances between the points
within the $3D$ objects in $\mathbb{R}^{3}$ and line integrals within
complex planes. We highlight the application of these concepts to
virtual tourism. 
\end{abstract}

\keywords{\textbf{Keywords: Riemannian metric, differential geometry, conformal
mapping, probability density functions, angle preservations, virtual
tourism, complex analysis.}}

\subjclass[2000]{MSC: 53B12, 30C20.}
\begin{center}
\vspace{1cm}
\par\end{center}

This article is part of:
\begin{center}
\vspace{0.1cm}
\par\end{center}

\textbf{\emph{Information Geometry, Volume 45: Handbook of Statistics,
Elsevier/North-Holland, Amsterdam }}\textbf{(2021 Fall)}

\section{\textbf{Introduction}}

C.R. Rao introduced his famous metric \cite{CRRao1949} in 1949 for
measuring distances between probability densities arising from population
parameters. This was later called by others the Rao distance (see,
for example, \cite{AtkinsonSankhya1981,Rios-}). There are several
articles available for the technicalities of Rao distance (see for
example, \cite{Amari-book,Jimenez,Chaudhuri,Chen,Nielsen-Arxiv})
and its applications (see for example, \cite{Rao-Krantz-Virtual,Taylor,Maybank}).
An elementary exposition of the same appeared during his centenary
in \cite{Planstino-Significance}. Rao distances and other research
contributions of renowned statistician C.R. Rao were recollected by
those who celebrated his 100th birthday during 2020 (see for example,
\cite{Efron-Amari-Rubin-Arni-Cox,BLSPRao et al,BLSP_ParthPP}). A
selected list of Rao's contributions in R programs was also made available
during his centenary (\cite{VinodHD}).

Rao distances are constructed under the framework of a quadratic differential
metric, Riemannian metric, and differential manifolds over probability
density functions and the Fisher information matrix. C.R. Rao considered
populations as abstract spaces which he called population spaces \cite{CRRao1949},
and then he endeavored to obtain topological distances between two
populations.

In the next section, we will describe manifolds. Section 3 will highlight
technicalities of Rao distances and Section 4 will treat conformal
mappings and basic constructions. Section 5 will conclude the chapter
with applications in virtual tourism.

\section{\textbf{Manifolds}}

Let $\mathbf{Df(a)}$ denotes the derivative of \textbf{$\mathbf{f}$
}at $\mathbf{a}$ for $\mathbf{a}\in\mathbb{R}^{n}$ and $\mathbf{f}:\mathbb{R}^{n}\rightarrow\mathbb{R}^{m}.$
A function $\mathbf{f}:\mathbb{R}^{n}\rightarrow\mathbb{R}^{m}$ is
\textit{differentiable} at $\mathbf{a}\in\mathbb{R}^{n}$ if there
exists a linear transformation $\mathbf{J}:\mathbb{R}^{n}\rightarrow\mathbb{R}^{m}$
such that

\begin{equation}
\lim_{h\rightarrow0}\frac{\left\Vert \mathbf{f(a}+\mathbf{h})-\mathbf{f(a)}-\mathbf{J}(\mathbf{h})\right\Vert }{\left\Vert \mathbf{h}\right\Vert }=0.\label{eq:1}
\end{equation}

Here $\mathbf{h}\in\mathbb{R}^{n}$ and $\mathbf{f(a+h)-f(a)-J(h)}\in\mathbb{R}^{n}.$
If $\mathbf{f}:\mathbb{R}^{n}\rightarrow\mathbb{R}^{m}$ is differentiable
at $a,$ then there exists a unique linear transformation $\mathbf{J}:\mathbb{R}^{n}\rightarrow\mathbb{R}^{m}$
such that (\ref{eq:1}) holds. The $m\times n$ matrix created by
$\mathbf{Df(a)}:\mathbb{R}^{n}\rightarrow\mathbb{R}^{m}$ is the Jacobian
matrix, whose elements are

\[
\begin{array}{cc}
\mathbf{Df(a)}= & \left[\begin{array}{cccc}
D_{1}f_{1}(a) & D_{2}f_{1}(a) & \cdots & D_{n}f_{1}(a)\\
D_{1}f_{2}(a) & D_{2}f_{2}(a) & \cdots & D_{n}f_{2}(a)\\
\vdots & \vdots &  & \vdots\\
D_{1}f_{m}(a) & D_{2}f_{m}(a) & \cdots & D_{n}f_{m}(a)
\end{array}\right]\end{array}
\]

That is, $\mathbf{J}(h)=\mathbf{Df(a)}.$ Since ${\bf J}$ is linear,
we have $\mathbf{J}(b_{1}\lambda_{1}+b_{2}\lambda_{2})=b_{1}\mathbf{J}(\lambda_{1})+b_{2}\mathbf{J}(\lambda_{2})$
for every $\mathbf{\lambda_{1},\lambda_{2}\in}\mathbb{R}^{n}$ and
every pair of scalars $b_{1}$ and $b_{2}.$ Also, the directional
derivative of \textbf{$\mathbf{f}$} at $\mathbf{a}$ in the direction
of $\mathbf{v}$ for $\mathbf{v}\in\mathbb{R}^{n}$ is denoted by
$\mathbf{D(f,v)}$ is given by

\begin{equation}
\mathbf{D(f,v)}=\lim_{h\rightarrow0}\frac{\left\Vert \mathbf{f(a+h{\bf v})}-\mathbf{f(a)}\right\Vert }{\left\Vert \mathbf{h}\right\Vert }\label{eq:2}
\end{equation}

\noindent provided \emph{R.H.S. }of (\ref{eq:2}) exists. When $\mathbf{f}$
is linear $\mathbf{D(f,v)=f(v)}$ for every $v$ and every $\mathbf{a}.$
Since $J(\mathbf{h})$ is linear, we can write

\begin{equation}
\mathbf{f(a+u)=f(a)+D(f,u)+\left\Vert u\right\Vert \Delta_{a}(u)},\label{eq:3}
\end{equation}

where

\begin{align*}
\mathbf{u} & \in\mathbb{R}^{n}\text{ with}\left\Vert \mathbf{u}\right\Vert <r\text{ for }r>0,\text{so that }\mathbf{a+u}\in\mathbf{B}(a;r)\\
 & \text{ \ensuremath{\qquad} for an \ensuremath{n-}ball }\mathbf{B}(a;r)\in\mathbb{R}^{n},
\end{align*}

\[
\mathbf{\Delta_{a}(u)}=\frac{\left\Vert \mathbf{f\mathbf{(a+h)-f(a)}}\right\Vert }{\left\Vert \mathbf{h}\right\Vert }-\mathbf{f'(a)}\text{ if }\mathbf{h}\neq0,
\]

\[
\mathbf{\Delta_{a}(u)}\rightarrow0\text{ as }\mathbf{u}\rightarrow0.
\]

When $\mathbf{u}=h\mathbf{v}$ in (\ref{eq:3}), we have

\begin{equation}
\mathbf{f(a+hv)-f(a)=hD(f,u)+}\left\Vert h\right\Vert \left\Vert \mathbf{v}\right\Vert \mathbf{\Delta_{a}(u)}\label{eq:4}
\end{equation}

For further results on the Jacobian matrix and differentiability properties,
refer to \cite{Tu-manifolds,Spivak-manifolds,Apostol-book}.

Consider a function $f=u+iv$ defined on the plane $\mathbb{C}$ with
$u(z),$ $v(z)\in\mathbb{R}$ for $z=(x,y)\in\mathbb{C}.$ If there
exists four partial derivatives

\begin{equation}
\frac{\partial u(x,y)}{\partial x},\frac{\partial v(x,y)}{\partial x},\frac{\partial u(x,y)}{\partial y},\frac{\partial v(x,y)}{\partial y},\label{eq:four partial derivatives}
\end{equation}

\noindent and these partial derivatives satisfy \emph{Cauchy-Riemann
}equations (\ref{eq:Cauchy-Riemann equations})

\begin{equation}
\frac{\partial u(x,y)}{\partial x}=\frac{\partial v(x,y)}{\partial y}\text{ and }\frac{\partial v(x,y)}{\partial x}=-\frac{\partial u(x,y)}{\partial y},\label{eq:Cauchy-Riemann equations}
\end{equation}

then

\[
\mathbf{D}f\mathbf{(a)}=\frac{\partial u(x,y)}{\partial x}+i\frac{\partial v(x,y)}{\partial x}\text{ for }u,v\in\mathbb{R}.
\]

\begin{thm}
Let $f=u(x,y)+iv(x,y)$ for $u(x,y),v(x,y)$ defined on a subset $B_{\mathbf{\delta}}(\mathbf{c)}$
$\subset\mathbb{C}$ for $\mathbf{\delta,c},(x,y)\in\mathbb{C}.$
Assume $u(x,y)$ and $v(x,y)$ are differentiable at an interior $\mathbf{a}=(a_{1},a_{2})\subset B_{\mathbf{\delta}}(\mathbf{c)}.$
Suppose the partial derivatives $\lim_{(x,y)\rightarrow\mathbf{a}}\frac{u(x,y)-u(\mathbf{a)}}{(x,y)-\mathbf{a}}$
and $\lim_{(x,y)\rightarrow\mathbf{a}}\frac{v(x,y)-v(\mathbf{a)}}{(x,y)-\mathbf{a}}$
exists for $\mathbf{a}$ and these partial derivatives satisfy Cauchy-Riemann
equations at $\mathbf{a.}$ Then

\[
\mathbf{D}f\mathbf{(a)}=\lim_{\mathbf{(u,v)\rightarrow(a_{1},a_{2})}}\frac{f(u,v)-f(\mathbf{a)}}{(u,v)-\mathbf{a}}
\]
exists, and

\[
\mathbf{D}f\mathbf{(a)}=\lim_{(x,y)\rightarrow\mathbf{a}}\frac{u(x,y)-u(\mathbf{a})}{(x,y)-\mathbf{a}}+i\left[\lim_{(x,y)\rightarrow\mathbf{a}}\frac{v(x,y)-v(\mathbf{a})}{(x,y)-\mathbf{a}}\right].
\]
\end{thm}

If $\mathbf{D}f\mathbf{(a)}$ exists for every $B_{\mathbf{\delta}}(\mathbf{c)}$
$\subset\mathbb{C}$ then we say that $f$ is holomorphic in $B_{\mathbf{\delta}}(\mathbf{c)}$
and is denoted as $H(B_{\mathbf{\delta}}(\mathbf{c)}).$ Readers are
reminded that when $f$ is a complex function in $B_{\mathbf{\delta}}(\mathbf{c)}$
$\subset\mathbb{C}$ that has a differential at every point of $B_{\mathbf{\delta}}(\mathbf{c)}$,
then $f\in H(B_{\mathbf{\delta}}(\mathbf{c)})$ if, and only if, the
Cauchy-Riemann equations (\ref{eq:Cauchy-Riemann equations}) are
satisfied for every $\mathbf{a}\in B_{\mathbf{\delta}}(\mathbf{c)}.$
Refer to \cite{Krantz-CAGeometric,Krantz-guide,Rudin-RealComplex,Apostol-book,Krantz-MATLAB}
for other properties of holomorphic functions and their association
with Cauchy-Riemann equations.

\subsection{Conformality between two regions}

Holomorphic functions discussed above allows us to study conformal
equivalences (i.e. angle preservation properties). Consider two regions
$B_{\mathbf{\delta}}(\mathbf{c)},B_{\alpha}(\mathbf{d)}\subset\mathbb{C}$
for some $\mathbf{c,d,}\delta,\alpha\in\mathbb{C}.$ These two regions
are conformally equivalence if there exists a function $g\in H(B_{\mathbf{\delta}}(\mathbf{c)})$
such that $g$ is one-to-one in $B_{\mathbf{\delta}}(\mathbf{c)}$
and such that $g(B_{\mathbf{\delta}}(\mathbf{c)})=B_{\alpha}(\mathbf{d)}.$
This means $g$ is conformally one-to-one mapping if $B_{\mathbf{\delta}}(\mathbf{c)}$
onto $B_{\alpha}(\mathbf{d)}$. The inverse of $g$ is holomorphic
in $B_{\alpha}(\mathbf{d)}$.

This implies $g$ is a conformal mapping of $B_{\alpha}(\mathbf{d)}$
onto $B_{\mathbf{\delta}}(\mathbf{c)}.$ We will introduce conformal
mappings in the next section. The two regions $B_{\mathbf{\delta}}(\mathbf{c)}$
and $B_{\alpha}(\mathbf{d)}$ are homeomorphic under the conformality.

The idea of manifolds is more general than the concept of a complex
plane. It uses the concepts of the Jacobian matrix, diffeomorphism
between $\mathbb{R}^{m}$ and $\mathbb{R}^{n}$, and linear transformations.
A set $M\subset\mathbb{R}^{n}$ is called a manifold if for every
$a\in M$, there exists a neighborhood $\mathbf{U}$ (open set) containing
$\mathbf{a}$ and a diffeomorphism $f_{1}:\mathbf{U}\rightarrow\mathbf{V}$
for $V\subset\mathbb{R}^{n}$ such that

\begin{equation}
\mathbf{f_{1}}(\mathbf{U}\cap M)=\mathbf{V}\cap\left(\mathbb{R}^{k}\times\{0\}\right)\label{eq:manifoldf1}
\end{equation}

The dimension of $M$ is $k.$ See \cite{Spivak-manifolds,Tu-manifolds}
for other details on manifolds. Further for an open set $\mathbf{V_{1}}\subset\mathbb{R}^{k}$
and a diffeomorphism 
\begin{equation}
\mathbf{f_{2}:V_{1}}\rightarrow\mathbb{R}^{n}\label{eq:f2}
\end{equation}

\noindent such that $\mathbf{Df_{2}(b)}$ has rank $k$ for $\mathbf{b\in}\mathbf{V}_{\mathbf{1}}.$ 
\begin{rem}
There exists a diffeomorphism as in (\ref{eq:f2}) such that $\mathbf{f_{2}}:\mathbf{V_{1}}\rightarrow\mathbf{f(V_{1})}$
is continuous. 
\end{rem}

\section{\textbf{Rao distance}}

A Riemannian metric is defined using an inner product function, manifolds,
and the tangent space of the manifold considered. 
\begin{defn}
\textbf{Riemannian metric: }Let $a\in M$ and $T_{a}M$ be the tangent
space of $M$ for each $a.$ A Riemannian metric $\mathcal{G}$ on
$M$ is an inner product

\[
\mathcal{G}_{a}:T_{a}M\times T_{a}M\rightarrow\mathbb{R}^{n}
\]

\noindent constructed on each $a.$ Here $(M,\mathcal{G})$ forms
Riemannian space or Riemannian manifold. The tensor space can be imagined
as collection of all the multilinear mappings from the elements in
$M$ as shown in Figure 3.1. For general references on metric spaces
refer to \cite{Kobayashi-Nomizu-book,Ambrisio-Luigi-book}.

\noindent 
\begin{figure}
\includegraphics{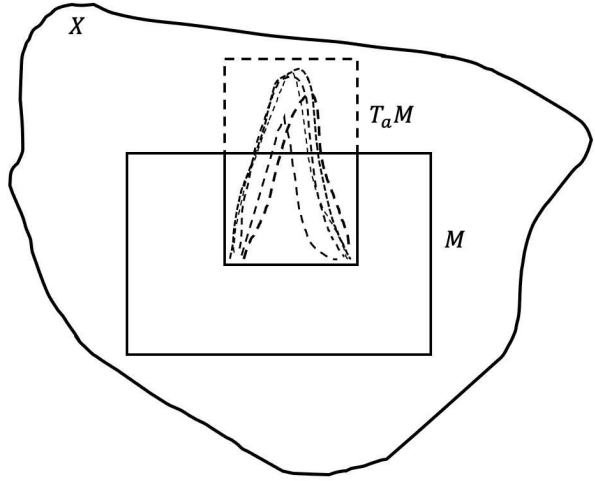}

\caption{\label{fig:MultilinearMappinginM}Mapping of elements in the manifold
$M$ in a metric space $X$ to the tensor space $T_{a}M.$}
\end{figure}
\end{defn}

Let $p(x,\theta_{1},\theta_{2},...,\theta_{n})$ be the probability
density function of a random variable $X$ such that $x\in X,$ and
$\theta_{1},\theta_{2},...,\theta_{n}$ are the parameters describing
the population. For different values of $\theta_{1},\theta_{2},...,\theta_{n}$
we will obtain different populations. Let us call $P(x,\Theta_{n})$
the population space created by $\Theta_{n}$ for a chosen functional
form of $X.$ Here $\Theta_{n}=\left\{ \theta_{1},\theta_{2},...,\theta_{n}\right\} .$
Let us consider another population space $P\left(x,\Theta_{n}+\Delta\right)$,
where 
\[
\Theta_{n}+\Delta=\left\{ \theta_{1}+\delta\theta_{1},\theta_{2}+\delta\theta_{2},...,\theta_{n}+\delta\theta_{n}\right\} .
\]

Let $\phi\left(x,\Theta_{n}\right)dx$ be the probability differential
corresponding to $P(x,\Theta_{n})$ and $\phi\left(x,\Theta_{n}+\Delta\right)dx$
be the probability differential corresponding to $P\left(x,\Theta_{n}+\Delta\right).$
Let

\begin{equation}
d\phi\left(\Theta_{n}\right)\label{eq:dphithea}
\end{equation}
be the differences in probability densities corresponding to $\Theta_{n}$
and $\Theta_{n}+\Delta.$ In (\ref{eq:dphithea}), C.R. Rao considered
only the first order differentials \cite{CRRao1949,Burbea-Rao,Micheccli-Noakes}.
The variance of the distribution of $\frac{d\phi}{\phi}$ is given
by

\begin{equation}
d\left[\frac{d\phi}{\phi}\right]^{2}=\sum\sum F_{ij}d\theta_{i}d\theta_{j}\label{eq:rao-metric}
\end{equation}

where $F_{ij}$ is the \emph{Fisher information matrix} for

\[
F_{ij}=E\left[\left(\frac{1}{\phi}\frac{\partial\phi}{\partial\theta_{i}}\right)\left(\frac{1}{\phi}\frac{\partial\phi}{\partial\theta_{j}}\right)\right]\text{ (for }E\text{ the expectation)}.
\]

Constructions in (\ref{eq:rao-metric}) and other measures between
probability distributions by C.R. Rao has played an important role
in statistical inferences.

Let $f_{3}$ be a measurable function on $X$ with differential $\phi\left(x,\Theta_{n}\right)dx.$
This implies that $f_{3}$ is defined on an interval $S\subset\mathbb{R}$
and there exists a sequence of step-functions $\{s_{n}\}$ on $S$
such that 
\[
\lim_{n\rightarrow\infty}s_{n}(x)=f_{3}(x)\text{ almost everywhere on }S
\]

for $x\in X.$

If $f_{3}$ is a $\sigma$-finite measure on $X$, then it satisfies

\[
\frac{d}{d\theta_{i}}\int_{S}P\left(x,\Theta_{n}\right)d\mu=\int_{S}\frac{dP\left(x,\Theta_{n}\right)}{d\theta}d\mu
\]

and

\[
\frac{d}{d\theta_{i}}\int_{S}P\left(x,\Theta_{n}\right)d\mu=\frac{d}{d\theta_{i}}\int_{S}\frac{P'\left(x,\Theta_{n}\right)}{P\left(x,\Theta_{n}\right)}P\left(x,\Theta_{n}\right)d\mu.
\]

\begin{rem}
Since the random variable $X$ can be covered by the collection of
sets $T_{n}$ such that

\[
\bigcup_{n=1}^{\infty}T_{n}=X,
\]
$\mu$ is the $\sigma$-finite measure, and

\[
f_{2}(x)>0\text{ and }\int f_{2}(x)\mu(dx)<\infty.
\]

\begin{figure}
\includegraphics{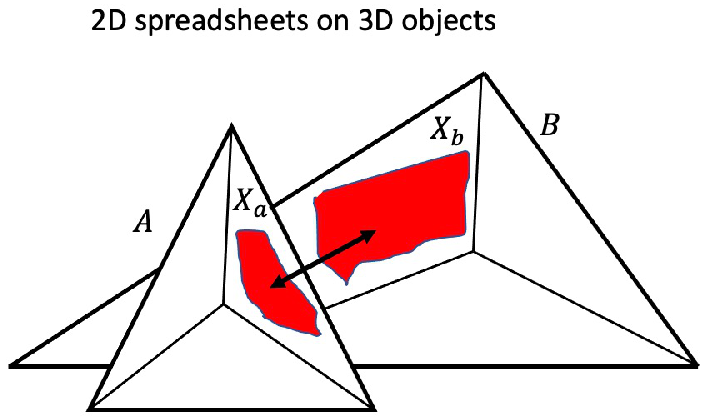}

\caption{\label{fig:Two-2D-shaped}Two 2D-shaped spreadsheets on 3D objects.
Metrics between such 2D-shaped spreadsheets can be studied based on
Rao distances. The distance between the space of points of $X_{a}$
located on the 3D shape $A$ to the space of $X_{b}$ located on the
3D object $B$ can be measured using population spaces conceptualized
in Rao distance.}
\end{figure}
\end{rem}

The idea of Rao distance can be used to compute the geodesic distances
between two 2D spreadsheets on two different 3D objects as shown in
Figure \ref{fig:Two-2D-shaped}. Burbea-Rao studied Rao distances
and developed $\alpha$- order entropy metrics for $\alpha\in\mathbb{R}$
\cite{Burbea-Rao}, given as

\begin{equation}
d\left[\frac{d\phi}{\phi}\right]_{\alpha}^{2}(\theta)=\sum_{i,j}^{n}=\mathcal{G}_{ij}^{(\alpha)}d\theta_{i}d\theta_{j}\label{eq:Burbea-Raometric}
\end{equation}

where

\begin{equation}
\mathcal{G}_{ij}^{(\alpha)}=\int_{X}P(x,\Theta_{n})^{\alpha}\left(\partial_{\theta_{i}}\log P\right)\left(\partial_{\theta_{j}}\log P\right)d\mu.\label{eq:Burbea-Rao-G}
\end{equation}

For the case of $P(x,\Theta_{n})$ as a multinomial distribution where
$x\in X$ for a sample space $X=\{1,2,...,n\},$ Burbea-Rao \cite{Burbea-Rao}
showed that

\begin{equation}
\mathcal{G}_{ij}^{(\alpha)}(\theta)=\int_{X}P(x,\Theta_{n})^{\alpha-2}\left(\partial_{\theta_{i}}\log P\right)\left(\partial_{\theta_{j}}\log P\right)d\mu.\label{eq:multinomial-tensor}
\end{equation}

The tensor of the metric in (\ref{eq:multinomial-tensor}) is of rank
$n.$

\section{\textbf{Conformal Mapping}}

\begin{landscape}

\begin{figure}
\includegraphics{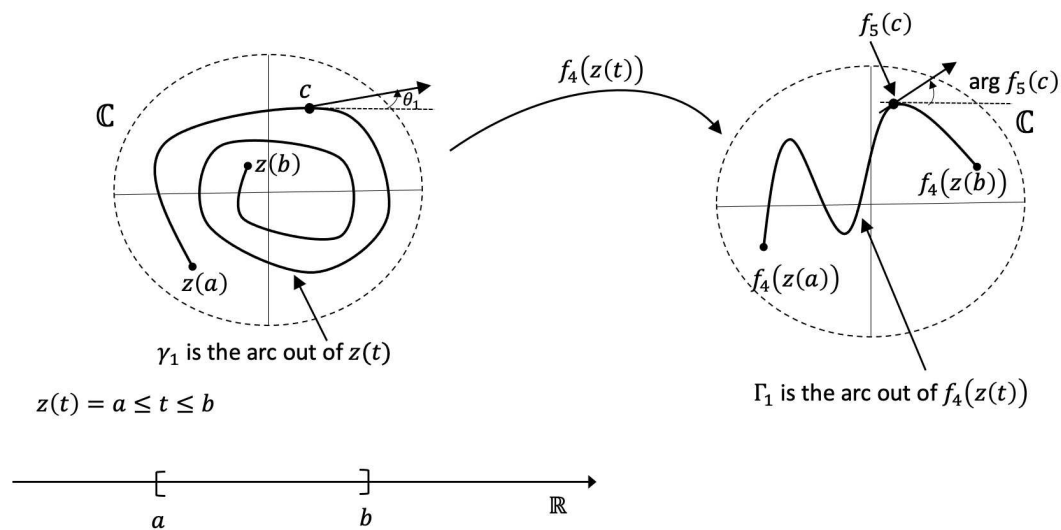}

\caption{\label{fig:Mapping-of-pointsbasics}Mapping of points from the real
line to an arc in the complex plane. Suppose $\gamma_{1}$ is the
\emph{arc} constructed out of $z(t)$ values. An arbitrary point $z(c)$
on $\gamma_{1}$ for $a\protect\leq c\protect\leq b$ at which $f_{4}$
is holomorphic and $f_{4}^{'}\left(z(c)\right)\protect\neq0.$ Let
$\theta_{1}$ be the angle of inclination at $c.$ When we denote
$\text{arg }f_{4}^{'}\left(z(c)\right)=\alpha_{1},$ it will lead
to $\text{arg }f_{5}^{'}(c)=\alpha_{1}+\theta_{1}$}
\end{figure}

\end{landscape}

\begin{landscape}
\begin{figure}
\includegraphics{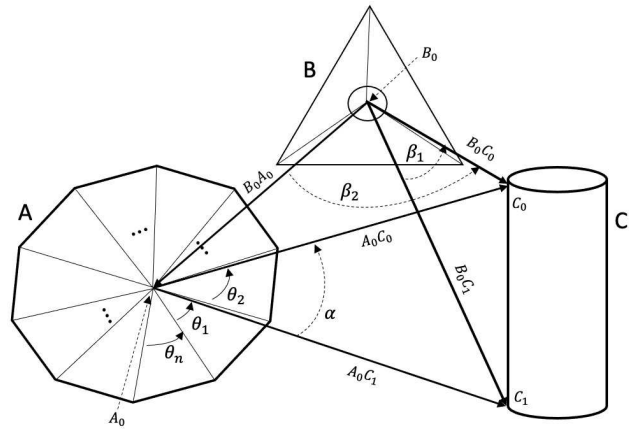}

\caption{\label{fig:3D-objects-and}3D objects and conformality with respect
to different viewpoints. The angles $\theta_{1}$, $\theta_{2}$,...
$\alpha$, $\beta_{1}$, $\beta_{2}$ are all measured. The distances
of the rays $A_{0}C_{0}$, $A_{0}C_{1}$, $B_{0}A_{0}$, $B_{0}C_{0}$,
$B_{0}C_{1}$ by assuming they are situated in a single $\mathbb{R}^{3}$
structure and also assuming they are situated in five different complex
planes is computed. By visualizing the three objects are replicas
of an actual tourist spot an application to virtual tourism is discussed
in section 5.}
\end{figure}

\end{landscape}

The storyline of this section is constructed around Figure \ref{fig:Mapping-of-pointsbasics}
and Figure \ref{fig:3D-objects-and}. First let us consider Figure
\ref{fig:Mapping-of-pointsbasics} for our understanding of conformal
mapping property. Let $z(t)$ be a complex-valued function for $z(t)=a\leq t\leq b$
for $a,b\in\mathbb{R}.$ Suppose $\gamma_{1}$ is the \emph{arc} constructed
out of $z(t)$ values. Suppose an arc $\Gamma_{1}$ is formed by the
mapping $f_{4}$ with a representation

\[
f_{5}(t)=f_{4}\left(z(t)\right)\text{ for }a\leq t\leq b.
\]

Let us consider an arbitrary point $z(c)$ on $\gamma_{1}$ for $a\leq c\leq b$
at which $f_{4}$ is holomorphic and $f_{4}^{'}\left(z(c)\right)\neq0.$
Let $\theta_{1}$ be the angle of inclination at $c$ as shown in
Figure \ref{fig:Mapping-of-pointsbasics}, then we can write $\text{arg }z^{'}(c)=\theta_{1}.$
Let $\alpha_{1}$ be the angle at $f_{4}^{'}\left(z(c)\right)$, i.e.

\[
\text{arg }f_{4}^{'}\left(z(c)\right)=\alpha_{1}.
\]

By this construction,

\begin{equation}
\text{arg }f_{5}^{'}(c)=\alpha_{1}+\theta_{1},\text{ (because arg }f_{5}^{'}(c)=\text{arg }f_{4}^{'}\left(z(c)\right)+\text{ arg }z^{'}(c))\label{eq:arg f5c}
\end{equation}

\noindent where $\text{arg }f_{5}^{'}(c)$ is the angle at $f_{5}^{'}(c)$
corresponding $\Gamma_{1}.$ Suppose that $\gamma_{2}$ is another
{arc} passing through $z(c)$ and $\theta_{2}$ be the angle of
inclination of the directed tangent line at $\gamma_{2}.$ Let $\Gamma_{2}$
be the {arc} corresponding to $\gamma_{2}$ and $\text{arg }f_{6}^{'}(c)$
be the corresponding angle at $f_{6}^{'}(c).$ Hence the two directed
angles created corresponding to $\Gamma_{1}$ and $\Gamma_{2}$ are

\begin{align*}
\text{arg }f_{5}^{'}(c) & =\alpha_{1}+\theta_{1}\\
\text{arg }f_{6}^{'}(c) & =\alpha_{2}+\theta_{2}
\end{align*}

This implies that

\begin{equation}
\text{arg }f_{6}^{'}(c)-\text{arg }f_{5}^{'}(c)=\theta_{2}-\theta_{1}.\label{eq:anglespreserved}
\end{equation}

The angle created from $\Gamma_{2}$ to $\Gamma_{1}$ at $f_{4}(z(c))$
is the same as the angle created at $c$ on $z(t)$ due to passing
of two {arcs} $\gamma_{1}$ and $\gamma_{2}$ at $c.$

Let $A,B,$ and $C$ be three $3D$ objects as shown in Figure \ref{fig:3D-objects-and}.
Object $A$ has a polygon-shaped structure with a pointed top located
at $A_{0}.$ A pyramid-shaped structure $B$ is located near object
$A$ and a cylinder-shaped object $C.$ Object $B$ has a pointed
top located at $B_{0}.$ Let $C_{0}$ be the nearest distance on $C$
from $B_{0}$ and \textbf{$C_{1}$ }be the farthest distance $C$
from $B_{0}.$ The norms of $A_{0},$ $B_{0},$ $C_{0}$, $C_{1}$
are all assumed to be different. Suppose $A_{0}=(A_{01,}A_{02},A_{03}),$
$B_{0}=(B_{01},B_{02},B_{03}),$ $C_{0}=(C_{01},C_{02},C_{03})$,
$C_{1}=(C_{11,}C_{12},C_{13}).$ Various distances between these points
are defined as below:

\begin{align}
A_{0}C_{0} & =\left\Vert A_{0}-C_{0}\right\Vert =\left[\sum_{i=1}^{3}\left(A_{0i}-C_{0i}\right)^{2}\right]^{1/2}\nonumber \\
A_{0}C_{1} & =\left\Vert A_{0}-C_{1}\right\Vert =\left[\sum_{i=1}^{3}\left(A_{0i}-C_{1i}\right)^{2}\right]^{1/2}\nonumber \\
B_{0}A_{0} & =\left\Vert B_{0}-A_{0}\right\Vert =\left[\sum_{i=1}^{3}\left(B_{0i}-A_{0i}\right)^{2}\right]^{1/2}\nonumber \\
B_{0}C_{0} & =\left\Vert B_{0}-C_{0}\right\Vert =\left[\sum_{i=1}^{3}\left(B_{0i}-C_{0i}\right)^{2}\right]^{1/2}\nonumber \\
B_{0}C_{1} & =\left\Vert B_{0}-C_{1}\right\Vert =\left[\sum_{i=1}^{3}\left(B_{0i}-C_{1i}\right)^{2}\right]^{1/2}\label{eq:fivedistances}
\end{align}

Let $\alpha$ be the angle from the ray $A_{0}C_{1}$ to the ray $A_{0}C_{0}$
with reference to the point $A_{0}$, $\beta_{1}$ be the angle from
the ray $B_{0}C_{1}$ to the ray $B_{0}C_{1}$ with reference to the
point $B_{0}$, and $\beta_{2}$be the angle from the ray $B_{0}A_{0}$
to the ray $B_{0}C_{0}$ with reference to the point $B_{0}$.

\begin{prop}
\label{prop:1}All the four points $A_{0},$ $B_{0},$ $C_{0}$, $C_{1}$
of Figure \ref{fig:3D-objects-and} can not be located in a single
Complex plane. These points could exist together in $\mathbb{R}^{3}.$ 
\end{prop}

\begin{proof}
Suppose the first coordinate of the plane represents the distance
from $x-$axis, the second coordinate is the distance from $y-$axis,
and the third coordinate represents the height of the 3D structures.
Even if $A_{03}=B_{03}=C_{03}$, still all the four points cannot
be on the same plane because $C_{03}$ cannot be equal to $C_{13}.$
Hence al the four points cannot be situated within a single complex
plane. However, by the same construction, they all can be situated
within a single 3D sphere or in $\mathbb{R}^{3}.$ 
\end{proof}
\begin{prop}
Suppose the norms and the third coordinates of $A_{0},$ $B_{0},$
$C_{0}$, $C_{1}$ are all assumed to be different. Then, it requires
five different complex planes, say, $\mathbb{C}_{1},$ $\mathbb{C}_{2}$,
$\mathbb{C}_{3}$, $\mathbb{C}_{4},$ and $\mathbb{C}_{5}$ such that
$A_{0},C_{0}\in\mathbb{C}_{1}$, $A_{0},C_{1}\in\mathbb{C}_{2}$,
$A_{0},B_{0}\in\mathbb{C}_{3}$, $B_{0},C_{0}\in\mathbb{C}_{4}$,
$B_{0},C_{1}\in\mathbb{C}_{5}.$ 
\end{prop}

\begin{proof}
By Proposition \ref{prop:1} all the four points $A_{0},$ $B_{0},$
$C_{0}$, $C_{1}$ cannot be in a single complex plane. Although the
third coordinates are different two out of four points can be considered
such that they fall within a same complex plane. Hence, the five rays
$A_{0}C_{0}$, $A_{0}C_{1}$, $B_{0}A_{0}$, $B_{0}C_{0}$, $B_{0}C_{1}$
can be accommodated in five different complex planes. 
\end{proof}
\begin{prop}
The angles $\alpha$, $\beta_{1}$, $\beta_{2}$ and five distances
of (\ref{eq:fivedistances}) are preserved when $A_{0},$ $B_{0},$
$C_{0}$, $C_{1}$ are situated together in $\mathbb{R}^{3}.$ 
\end{prop}

\begin{proof}
The angle $\alpha$ is created while viewing the $3D$ structure $C$
from point $A_{0}.$ The angle $\beta_{1}$ is created while viewing
the $3D$ structure $C$ from the point $B_{0}.$ The angle $\beta_{2}$
is created while viewing the $3D$ structure $C$ from the point $A_{0}.$
These structures could be imagined to stand on a disc within a $3D$
sphere or in $\mathbb{R}^{3}$ even proportionately mapped to $\mathbb{R}^{3}.$
Under such a construction, without altering the ratios of various
distances, the angles remain the same in the mapped $\mathbb{R}^{3}.$ 
\end{proof}
Let us construct an arc $A_{0}C_{0}(t_{1})=a_{1}\leq t_{1}\leq b_{1}$
from the point $A_{0}$ to $C_{0}$ and call this arc $C_{1}.$ Here
$a_{1},b_{1}\in\mathbb{R}$ and $A_{0},C_{0}\in\mathbb{C}_{1}$. The
points of $C_{1}$ are $A_{0}C_{0}(t_{1}).$ The values of $t_{1}$
can be generated using a parametric representation which could be
a continuous random variable or a deterministic model. 
\begin{equation}
t_{1}=\psi_{1}(\tau)\text{ for }\alpha_{1}\leq\tau\leq\beta_{1}.\label{eq:parametric_t1}
\end{equation}

Then the arc length $L(C_{1})$ for the arc $C_{1}$ is obtained through
the integral

\begin{equation}
L(C_{1})=\int_{\alpha_{1}}^{\beta_{1}}\left|A_{0}C_{0}^{'}\left[\psi_{1}(\tau)\right]\right|\psi_{1}^{'}(\tau)d\tau.\label{eq:L(c1)}
\end{equation}

Likewise, the arc lengths $L(C_{2}),$ $L(C_{3})$, $L(C_{4})$, $L(C_{5})$
for the arcs $C_{2}$, $C_{3}$, $C_{4}$, $C_{5}$ are constructed
as follows:

\begin{equation}
L(C_{2})=\int_{\alpha_{2}}^{\beta_{2}}\left|A_{0}C_{1}^{'}\left[\psi_{2}(\tau)\right]\right|\psi_{2}^{'}(\tau)d\tau,\label{eq:L(c2)}
\end{equation}

where $A_{0}C_{1}(t_{2})=a_{2}\leq t_{2}\leq b_{2}$ for $a_{2},b_{2}\in\mathbb{R}$
and $A_{0},C_{1}\in\mathbb{C}_{2}$ and with parametric representation
$t_{2}=\psi_{2}(\tau)\text{ for }\alpha_{2}\leq\tau\leq\beta_{2}.$

\begin{equation}
L(C_{3})=\int_{\alpha_{3}}^{\beta_{3}}\left|B_{0}A_{0}^{'}\left[\psi_{3}(\tau)\right]\right|\psi_{3}^{'}(\tau)d\tau,\label{eq:LC3}
\end{equation}

where $B_{0}A_{0}(t_{3})=a_{3}\leq t_{3}\leq b_{3}$ for $a_{3},b_{3}\in\mathbb{R}$
and $B_{0},A_{0}\in\mathbb{C}_{3}$ and with parametric representation
$t_{3}=\psi_{3}(\tau)\text{ for }\alpha_{3}\leq\tau\leq\beta_{3}.$

\begin{equation}
L(C_{4})=\int_{\alpha_{4}}^{\beta_{4}}\left|B_{0}C_{0}^{'}\left[\psi_{4}(\tau)\right]\right|\psi_{4}^{'}(\tau)d\tau,\label{eq:LC4}
\end{equation}

where $B_{0}C_{0}(t_{4})=a_{4}\leq t_{4}\leq b_{4}$ for $a_{4},b_{4}\in\mathbb{R}$
and $B_{0},C_{0}\in\mathbb{C}_{4}$ and with parametric representation
$t_{4}=\psi_{4}(\tau)\text{ for }\alpha_{4}\leq\tau\leq\beta_{4}.$

\begin{equation}
L(C_{5})=\int_{\alpha_{5}}^{\beta_{5}}\left|B_{0}C_{1}^{'}\left[\psi_{5}(\tau)\right]\right|\psi_{5}^{'}(\tau)d\tau,\label{eq:Lc5}
\end{equation}

where $B_{0}C_{1}(t_{5})=a_{5}\leq t_{5}\leq b_{5}$ for $a_{5},b_{5}\in\mathbb{R}$
and $B_{0},C_{1}\in\mathbb{C}_{5}$ and with parametric representation
$t_{5}=\psi_{5}(\tau)\text{ for }\alpha_{5}\leq\tau\leq\beta_{5}.$ 
\begin{rem}
One could also consider a common parametric representation

\[
\psi_{i}(\tau)=\psi(\tau)\text{ for }i=1,2,...,5
\]
if that provides more realistic situation of modeling. 
\end{rem}

\section{\textbf{Applications}}

The angle preservation approach can be used in preserving the angles
and depth of 3D images for actual 3D structures. Earlier Rao \& Krantz
\cite{Rao-Krantz-Virtual} proposed such measures in the virtual tourism
industry.

Advanced virtual tourism technology is in the early stage of development
and it occupies a small fraction of the total tourism-related business.
Due to the pandemics and other large-scale disruptions around tourist
locations, there will be a high demand for virtual tourism facilities.
One such was visualized during COVID-19 (\cite{Rao-Krantz-Virtual}).
Let us consider a tourist location that has three 3D structured buildings
as in Figure 4.2. When a tourist visits the location in person then
such scenery can be seen directly from the ground level by standing
in between the three structures or standing beside one of the structures.
It is not always possible to see those features when standing above
those buildings. Suppose a video recording is available that was recorded
with regular video cameras; then the distances $A_{0}C_{0}$, $A_{0}C_{1}$,
$B_{0}A_{0}$, $B_{0}C_{0}$, $B_{0}C_{1}$ and angles $\alpha,$
$\beta_{1}$, $\beta_{2}$ would not be possible to capture. That
depth of the scenery and relative elevations and distances would not
be accurately recorded. The in-person virtual experience at most can
see the distance between the bottom structures of the tourist attractions.

The same scenery of Figure 4.2, when watched in person at some time
of the day, would be different when it is watched at a different time
due to the differences between day and night visions. The climatic
conditions and weather would affect the in-person tourism experiences.
All these can be overcome by having virtual tourism technologies proposed
for this purpose \cite{Rao-Krantz-Virtual}. The new technology called
LAPO (live-streaming with actual proportionality of objects) would
combine the pre-captured videos and photos with live-streaming of
the current situations using advanced drone technology. This would
enhance the visual experience of live videos by mixing them with pre-recorded
videos. Such technologies will not only enhance the visualizations
but also help in repeated seeing of the experiences and a closer look
at selected parts of the videos. Mathematical formulations will assist
in maintaining the exactness and consistency of the experiences. We
hope that the newer mathematical constructions, theories, and models
will also emerge from these collaborations.

The line integrals $L(C_{i})$ for $i=1,2,...,5$ are computed and
the angles between the structures can be practically pre-computed
for each tourist location so that these can be mixed with the live
streaming of the tourist locations. The angle preservation capabilities
to maintain the angles between various base points can be preserved
with actual measurements that will bring a real-time experience of
watching the monuments.

The virtual tourism industry has many potential advantages if it is
supported by high-end technologies. Viewing the normal videos of tourist
attractions through the internet browser could be enriched with the
new technology proposed \cite{Rao-Krantz-Virtual}. These new technologies
combined with more accurate preservations of the depth, angles, and
relative distances would enhance the experiences of virtual tourists.
Figure 4.2 could be considered as a view of a tourist location. There
are more realistic graphical descriptions available to understand
the proposed technology LAPO using the information geometry and conformal
mapping \cite{Rao-Krantz-Virtual}.

Apart from applying mathematical tools, there are advantages of virtual
tourism. Although this discussion is out of scope for this article,
we wish to highlight below a list of advantages and disadvantages
of new virtual tourism technology taken from \cite{Rao-Krantz-Virtual}.
\medskip{}
\\

\textsc{Advantages:} 
\begin{enumerate}
\item[\textbf{(a)}] Environmental protection around ancient monuments; 
\item[\textbf{(b)}] Lesser disease spread at the high population density tourist locations; 
\item[\textbf{(c)}] Easy tour for physically challenged persons; 
\item[\textbf{(d)}] Creation of newer employment opportunities; 
\item[\textbf{(e)}] The safety of tourists; 
\item[\textbf{(f)}] The possibility of the emergence of new software technologies. 
\end{enumerate}
\vspace*{0.15in}

\textsc{Disadvantages:} 
\begin{enumerate}
\item[\textbf{(a)}] Possible abuse of the technology that can harm the environment around
the tourist locations; 
\item[\textbf{(b)}] Violation of individual privacy; 
\item[\textbf{(c)}] Misuse of drone technology. 
\end{enumerate}
\vspace*{0.15in}

Overall there are plenty of advantages of developing this new technology
and implementing it with proper care taken for protection against
misuse. The importance of this technology is that it will have deeper
mathematical principles and insights that were not utilized previously
in the tourism industry. When the population mobility reduces due
to pandemics the hospitality and business industry was seen to have
severe financial losses. In such a situation, virtual tourism could
provide an alternative source of financial activity.

There are of course several advantages of real tourism too, like understanding
the actual physical structures of the monuments, touching of the monuments
(trees, stones, water, etc.,), and feeling real climatic conditions.
We are not describing here all the possible advantages and disadvantages
between virtual versus real tourism experiences.

The concept of Rao distance constructed on population spaces can be
used to measure distances between two probability densities. One possible
application is to virtual tourism. This article is anticipated to
help understand various technicalities of Rao distances and conformal
mappings in a clear way.

\subsection*{Acknowledgements: }

ASRS Rao thanks to his friend Padala Ramu who taught him complex analysis
and to all the students who had attended ASRSR's courses on real and
complex analysis.

\end{document}